\documentclass[11pt]{article}
\usepackage[russian]{babel}
\usepackage[cp1251]{inputenc}
\usepackage{amsmath,latexsym}
\usepackage[psamsfonts]{amssymb}
\usepackage[dvips]{graphicx}
\usepackage[mathcal]{euscript}
\begin{document}

\begin{center}
\textbf{Initial boundary value problem for a high-order equation \\ with two lines of degeneracy with the Caputo derivative
 }\\
\textbf{B.Yu.Irgashev}\\
\emph{Namangan  Engineering Construction Institute, \\ Institute of Mathematics named after V.I.Romanovsky of the Academy of Sciences of the Republic of Uzbekistan.}\\
  \emph{E-mail: bahromirgasev@gmail.com}
\end{center}

\textbf{Abstract.} \emph{ In the article, in a rectangular domain, by the Fourier method, the initial boundary value problem for a high-order equation with two lines of degeneracy with a fractional derivative in the sense of Caputo is investigated for uniqueness and solvability. Sufficient conditions for the well-posedness of the problem posed are obtained.}

\\
\section*{1. Introduction and problem statement }

  Fractional partial differential equations underlie the mathematical modeling of various physical processes and environmental phenomena that have a fractal nature [1-2]. Model equations with fractional derivatives with constant coefficients are well studied. Recently, specialists have intensively studied equations with variable coefficients. Degenerate equations are among such equations. The number of works on
degenerate equations with fractional derivatives is relatively small. Many of these papers deal with ordinary differential equations. In
the work [2] the equation
\[D_{0x}^\alpha {t^\beta }u\left( t \right) = \lambda u\left( x \right),\,\,0 < x < b,\]
where $0 < \alpha < 1,\,\,\lambda $ is the  spectral parameter, $\beta = const \ge 0.$
In the same work, it was noted that the equation
\[D_{0x}^\alpha {t^\beta }u\left( t \right) + \sum\limits_{j = 1}^n {{a_j}\left( x \right)} D_{0x}^{{\alpha _j}}u\left( t \right) + b\left( x \right)u = c\left( x \right),0 \ne {\alpha _j} < \alpha ,\]
plays an important role in the theory of inverse problems for degenerate equations of hyperbolic type. In the article [3] the following problem was investigated for solvability
\[\left\{ \begin{array}{l}
\partial _0^\nu \left( {k\left( t \right)y\left( t \right)} \right) + c\left( {y\left( t \right)} \right) = f\left( t \right),\\
k\left( 0 \right)y\left( 0 \right) = 0,
\end{array} \right.\]
Here $\partial _0^\nu  - $ is fractional differentiation operator in the Caputo sense, $0 < \nu  < 1,$ $k\left( t \right) \in {C^1}\left[ {0,T} \right],\,\,\,k\left( t \right) \ge 0,$   for $t \in \left[ {0;T} \right],$  $c\left( \eta  \right) \in C\left( R \right),\,\,c\left( 0 \right) = 0.$
In the paper [4] solutions in the closed form of the fractional order equations were found
\[\left( {D_{0 + }^\alpha y} \right)\left( x \right) = a{x^\beta }y\left( x \right) + f\left( x \right)\left( {0 < x < d \le \infty ,\alpha  > 0,\beta  \in R,a \ne 0} \right),\]
\[\left( {D_ - ^\alpha y} \right)\left( x \right) = a{x^\beta }y\left( x \right) + f\left( x \right)\left( {0 \le d < x < \infty ,\alpha  > 0,\beta  \in R,a \ne 0} \right),\]
with fractional Riemann-Liouville derivatives on the semiaxis $\left( {0,\infty } \right)$ [5]:
\[\left( {D_{0 + }^\alpha y} \right)\left( x \right) = {\left( {\frac{d}{{dx}}} \right)^{[\alpha ] + 1}}\frac{1}{{\Gamma \left( {1 - \{ \alpha \} } \right)}}\int\limits_0^x {\frac{{y\left( t \right)dt}}{{{{\left( {x - t} \right)}^{\{ \alpha \} }}}}} ,\left( {x > 0;\alpha  > 0} \right),\]
\[\left( {D_ - ^\alpha y} \right)\left( x \right) = {\left( { - \frac{d}{{dx}}} \right)^{[\alpha ] + 1}}\frac{1}{{\Gamma \left( {1 - \{ \alpha \} } \right)}}\int\limits_x^\infty  {\frac{{y\left( t \right)dt}}{{{{\left( {t - x} \right)}^{\{ \alpha \} }}}}} ,\left( {x > 0;\alpha  > 0} \right),\]
($[\alpha ]$ and $\{ \alpha \} $ mean the integer and fractional parts of the real number $\alpha $ ). Applied problems lead to such equations [6]. An example of such an equation gives the equation of the theory of polarography [7]
\[\left( {D_{0 + }^{{1 \mathord{\left/
 {\vphantom {1 2}} \right.
 \kern-\nulldelimiterspace} 2}}y} \right)\left( x \right) = a{x^\beta }y\left( x \right) + {x^{ - {1 \mathord{\left/
 {\vphantom {1 2}} \right.
 \kern-\nulldelimiterspace} 2}}},\left( {0 < x, - {1 \mathord{\left/
 {\vphantom {1 {2 < \beta  \le }}} \right.
 \kern-\nulldelimiterspace} {2 < \beta  \le }}0} \right),\]
 which arises for $a = - 1$ in diffusion problems [7].
As for degenerate partial differential equations involving fractional derivatives, it should be noted that research in this area is quite new. We note here the works [8-17], in which initial and various boundary value problems for degenerate equations with fractional derivatives were studied.

In this paper, in a rectangular domain , we study a boundary value problem for an equation of high even order involving a fractional derivative , which has degeneracy in both variables.

In the region $\Omega = {\Omega _x} \times {\Omega _y},\,{\mkern 1mu} {\kern 1pt} {\Omega _x} = \left\{ {x:{\mkern 1mu} { \kern 1pt} {\mkern 1mu} {\kern 1pt} 0 < x < 1} \right\},{\mkern 1mu} {\kern 1pt} \,{\Omega _y} = \left\{ {y: {\mkern 1mu} {\kern 1pt} {\mkern 1mu} {\kern 1pt} 0 < y < 1} \right\},$ consider the equation
\[_CD_{0x}^\alpha u\left( {x,y} \right) + {x^\beta}K\left( y \right)l\left( {u\left( {x,y} \right)} \right) = 0,\eqno(1)\]
where
\[l\left( {u\left( {x,y} \right)} \right) = {\left( { - 1} \right)^s}\frac{{{\partial ^{2s}}u\left( {x,y} \right)}}{{\partial {y^{2s}}}} + \frac{{{\partial ^{s - 1}}}}{{\partial {y^{s - 1}}}}\left( {{{\left( { - 1} \right)}^{s - 1}}{p_{s - 1}}\left( y \right)\frac{{{\partial ^{s - 1}}u\left( {x,y} \right)}}{{\partial {y^{s - 1}}}}} \right) + ...\]
\[ + \frac{\partial }{{\partial y}}\left( { - {p_1}\left( y \right)\frac{{\partial u\left( {x,y} \right)}}{{\partial y}}} \right) + {p_0}\left( y \right)u\left( {x,y} \right),\]
\[ {p_j}\left( y \right) \in {C^j}\left( {{{\overline \Omega  }_y}} \right),{\mkern 1mu} {\kern 1pt} {\mkern 1mu} {\kern 1pt} j = 0,1,...,s - 1,\,s \in N,\]
$K(y) \in C[0;1]$, $K(y)$ a sufficiently smooth function for $y \in (0;1]$,
$$K\left( y \right) > 0,y \in \left( {0,1} \right],K(0)=0, {K^{\left( i \right)}}\left( y \right) = O\left( {{y^{m -i}}} \right),y \to  + 0,0 \le m < s,i=0,1,...,$$
$0 < \alpha  < 1, - \alpha  < \beta \in R,{}_CD_{0x}^\alpha  - $ fractional Caputo differentiation operator of order $\alpha $
\[{}_CD_{0x}^\alpha u\left( {x,y} \right) = \frac{1}{{\Gamma \left( {1 - \alpha } \right)}}\int\limits_0^x {\frac{{\frac{\partial }{{\partial \tau }}u\left( {\tau ,y} \right)d\tau }}{{{{\left( {x - \tau } \right)}^\alpha }}}.} \]
For equation (1) consider the problem.

\textbf{Problem A.} Find a solution to equation (1) from the class
\[{}_CD_{0x}^\alpha u\left( {x,y} \right) \in C\left( \Omega  \right),u\left( {x,y} \right) \in C\left( {\overline \Omega  } \right),\]
\[\frac{{{\partial ^{2s - 1}}u\left( {x,y} \right)}}{{\partial {y^{2s - 1}}}} \in C\left( {\overline \Omega  } \right),{\mkern 1mu} {\kern 1pt} \frac{{{\partial ^{2s}}u\left( {x,y} \right)}}{{\partial {y^{2s}}}} \in C\left( \Omega  \right),\]
satisfying the conditions
$$
\frac{{{\partial ^j}u\left( {x,0} \right)}}{{\partial {y^j}}} = \frac{{{\partial ^j}u\left( {x,1} \right)}}{{\partial {y^j}}} = 0,{\mkern 1mu} {\kern 1pt} {\mkern 1mu} {\kern 1pt} {\mkern 1mu} {\kern 1pt} 0 \le x \le 1,{\mkern 1mu} {\kern 1pt} j = 0,1,...,s - 1,\eqno(2)
$$
$$
u\left( {0,y} \right) = \varphi \left( y \right),\eqno(3)
$$
here the function $\varphi \left( y \right) - $ is sufficiently smooth and satisfies the matching conditions.
\section*{2. Finding a Solution }
We are looking for a solution in the form
\[u\left( {x,y} \right) = X\left( x \right)Y\left( y \right).\]
Then, with respect to the variable $y$ , taking into account condition (2), we obtain the following spectral problem:
\[\left\{ {\begin{array}{*{20}{l}}
{l\left( {Y\left( y \right)} \right) = \lambda \frac{{Y(y)}}{{K\left( y \right)}},}\\
{{Y^{\left( j \right)}}(0) = {Y^{\left( j \right)}}\left( 1 \right) = 0,j = 0,1,...,s - 1.}
\end{array}} \right.\eqno(4)\]
We reduce problem (4) to an integral equation using the Green's function and obtain the necessary estimates for the eigenfunctions. Further, we will assume that $\lambda > 0$, this condition will be satisfied, for example, in the case when $0 \le {p_j}\left( y \right),j=0,1,...,s-1.$ This implies the existence of a continuous , symmetric Green's function, which will be denoted by $G\left( {y,\xi } \right).$ Taking into account the boundary conditions at the point $y=0$ and applying the Lagrange theorem on finite increments, we have
\[{Y}(y) = O\left( {{y^{s}}} \right),{\mkern 1mu} {\kern 1pt} y \to  + 0,\]
this relation is also valid for the Green's function. It remains to show the existence of eigenvalues and eigenfunctions of problem (4). The integral equation equivalent to problem (4) has the form
\[Y\left( y \right) = \lambda \int\limits_0^1 {\frac{{G\left( {y,\xi } \right)Y\left( \xi  \right)d\xi }}{{K\left( \xi  \right)}}} ,\eqno(5)\]
we write (5) in the form
\[\frac{{Y\left( y \right)}}{{\sqrt {K\left( y \right)} }} = \lambda \int\limits_0^1 {\frac{{G\left( {y,\xi } \right)}}{{\sqrt {K\left( \xi  \right)} \sqrt {K\left( y \right)} }}\frac{{Y\left( \xi  \right)}}{{\sqrt {K\left( \xi  \right)} }}d\xi } ,\]
we introduce the notation
\[\overline Y \left( y \right) = \frac{{Y\left( y \right)}}{{\sqrt {K\left( y \right)} }},\overline G \left( {y,\xi } \right) = \frac{{G\left( {y,\xi } \right)}}{{\sqrt {K\left( \xi  \right)} \sqrt {K\left( y \right)} }},\]
from this we have
$$
\bar Y\left( y \right) = \lambda \int\limits_0^1 {\bar G\left( {y,\xi } \right)\bar Y\left( \xi  \right)d\xi }.\eqno(6)$$
(6) is an integral equation with a symmetric and continuous, in both variables, kernel. According to the theory of equations with symmetric kernels, equation (6) has at most a countable number of eigenvalues and eigenfunctions. Thus, problem (4) has eigenvalues ${\lambda _n} > 0,n = 1,2,...,$ and the corresponding eigenfunctions are ${Y_n}\left( y \right).$ Let's place the eigenvalues in the following order $0 < {\lambda _1} < {\lambda _2} < ....$ Further, we will assume that
\[{\left\| {{Y_n}(y)} \right\|^2} = \int\limits_0^1 {\frac{{Y_n^2(y)dy}}{{K\left( y \right)}}}  = 1,\]
from here, taking into account (5), we have the Bessel inequality
\[\sum\limits_{n =1}^\infty  {{{\left( {\frac{{{Y_n}\left( y \right)}}{{{\lambda _n}}}} \right)}^2}}  \le \int\limits_0^1 {\frac{{{G^2}\left( {y,\xi } \right)d\xi }}{{K\left( \xi  \right)}}}  < \infty.\eqno(7)\]
Let us find the conditions under which the initial function is expanded into a Fourier series in terms of the system of eigenfunctions of problem (4). The following theorem holds.

\textbf{Theorem 1.} If the function $\varphi \left( y \right)$ satisfies the following conditions:

$1).  {\varphi ^{\left( j \right)}}\left( 0 \right) = {\varphi ^{\left( j \right)}}\left( 1 \right) = 0,j = \overline {0,s - 1};$

2). ${\varphi ^{\left( {2s} \right)}}\left( y \right)$ is continuous on $\left[ {0,1} \right] ,$

then it can be expanded in terms of the eigenfunctions of the problem (4), which will converge uniformly and absolutely.

\textbf{ Proof.} Let $l\left( {\varphi \left( y \right)} \right) = g\left( y \right),$ then
\[\varphi \left( y \right) = \int\limits_0^1 {G\left( {y,\xi } \right)g\left( \xi  \right)} d\xi ,\]
from here
\[\frac{{\varphi \left( y \right)}}{{\sqrt {K\left( y \right)} }} = \int\limits_0^1 {\frac{{G\left( {y,\xi } \right)}}{{\sqrt {K\left( y \right)K\left( \xi  \right)} }}} \sqrt {K\left( \xi  \right)} g\left( \xi  \right)d\xi ,\]
or
\[\frac{{\varphi \left( y \right)}}{{\sqrt {K\left( y \right)} }} = \int\limits_0^1 {\overline G \left( {y,\xi } \right)} \sqrt {K\left( \xi  \right)} g\left( \xi  \right)d\xi.\]
Let us apply the Hilbert-Schmidt theorem, then we will have
\[\frac{{\varphi \left( y \right)}}{{\sqrt {K\left( y \right)} }} = \sum\limits_{n = 1}^\infty  {{c_n}\frac{{{Y_n}\left( y \right)}}{{\sqrt {K\left( y \right)} }}} ,\]
where
\[{c_n} = \int\limits_0^1 {\frac{{\varphi \left( y \right){Y_n}\left( y \right)}}{{K\left( y \right)}}} dy.\]
Further, after reduction by $\frac{1}{{\sqrt {K\left( y \right)} }}$, we get
\[\varphi \left( y \right) = \sum\limits_{n = 1}^\infty  {{c_n}{Y_n}\left( y \right)} .\]

\textbf{Theorem 1 is proved.}

We fix $n$ and start solving the problem of finding the function $X_{n}(x)$. We will assume that the initial function $\varphi \left( y \right)$ satisfies the conditions of Theorem 1. Taking into account condition (3), we obtain the following initial problem:
\[\left\{ {\begin{array}{*{20}{l}}
{{}_CD_{0x}^\alpha {X_n}\left( x \right) =  - {\lambda _n}{x^\beta }{X_n}\left( x \right),}\\
{{X_n}\left( 0 \right) = {\varphi _n},}
\end{array}} \right.\eqno(8)\]
where
\[{\varphi _n} = \int\limits_0^1 {\frac{{\varphi \left( y \right){Y_n}\left( y \right)dy}}{{K\left( y \right)}}} ,{\mkern 1mu} {\kern 1pt} n = 1,2,....\]
Using the results of [14], the solution to problem (8) is written as
\[{X_n}\left( x \right) = {\varphi _n}{E_{\alpha ,\frac{\beta }{\alpha } + 1,\frac{\beta }{\alpha }}}\left( { - {\lambda _n}{x^{\alpha  + \beta }}} \right),\]
where
\[{E_{\alpha ,m,l}}\left( z \right) = \sum\limits_{i = 0}^\infty  {{c_i}{z^i}} ,{c_0} = 1,{c_i} = \prod\limits_{j = 0}^{i - 1} {\frac{{\Gamma \left( {\alpha \left( {jm + l} \right) + 1} \right)}}{{\Gamma \left( {\alpha \left( {jm + l + 1} \right) + 1} \right)}}} ,i \ge 1\]
is a Kilbas-Saigo function  [14].\\
It follows from the results of [18] that the following estimate holds:
\[\left| {{E_{\alpha ,\frac{\beta }{\alpha } + 1,\frac{\beta }{\alpha }}}\left( { - {\lambda _n}{x^{\alpha  + \beta }}} \right)} \right| \le M = const.\eqno(9)\]
Now we will seek the solution of the problem in the form
$$
u\left( {x,y} \right) = \sum\limits_{n = 1}^\infty  {{X_n}\left( x \right){Y_n}\left( y \right)}.\eqno(10)
$$
Let us find the conditions under which (10) is a regular solution to Problem A. The following theorem is true.

\textbf{Theorem 2.} Let the function $\varphi \left( y \right),$ satisfy the following conditions:	

$1).   \varphi \left( y \right) \in {C^{2k}}\left[ {0,1} \right],{\mkern 1mu} {\kern 1pt} {\varphi ^{\left( j \right)}}\left( 0 \right) = {\varphi ^{\left( j \right)}}\left( 1 \right) = 0,j = 0,1,...,k - 1,$

$2).{\left( {K\left( y \right)l\left( {\varphi \left( y \right)} \right)} \right)^{\left( j \right)}}\left( 0 \right) = {\left( {K\left( y \right)l\left( {\varphi \left( y \right)} \right)} \right)^{\left( j \right)}}\left( 1 \right) = 0,K\left( y \right)l\left( {\varphi \left( y \right)} \right) \in {C^{2s}}\left[ {0,1} \right],{\mkern 1mu} {\kern 1pt} {\mkern 1mu} {\kern 1pt} j = 0,1,...,s - 1,$
	
then a solution to problem A exists.

\textbf{ Proof.} (10) formally satisfies equation (1). Let us show the possibility of differentiating the series. Considering (9) we have
\[\left| {u\left( {x,y} \right)} \right| \le M {\mkern 1mu} {\kern 1pt} \sum\limits_{n = 1}^\infty  {\left| {{\varphi _n}} \right|\left| {{Y_n}\left( y \right)} \right|} ,\]
we apply the Cauchy-Bunyakovsky inequality
\[\sum\limits_{n = 1}^\infty  {\left| {{\varphi _n}} \right|\left| {{Y_n}\left( y \right)} \right|}  \le \sqrt {\sum\limits_{n = 1}^\infty  {{{\left( {\frac{{{Y_n}\left( y \right)}}{{{\lambda _n}}}} \right)}^2}} } \sqrt {\sum\limits_{n = 1}^\infty  {{{\left( {{\lambda _n}{\varphi _n}} \right)}^2}} }.\]
The convergence of the first multiplier in the last inequality follows from (7) Consider the second multiplier. We have
\[{\varphi _n} = \int\limits_0^1 {\frac{{\varphi \left( y \right){Y_n}(y)dy}}{{K\left( y \right)}}}  = \frac{1}{{{\lambda _n}}}\int\limits_0^1 {\varphi \left( y \right)l\left( {{Y_n}\left( y \right)} \right)dy} .\]
Integrating by parts, we get
\[{\varphi _n} = \frac{1}{{{\lambda _n}}}\int\limits_0^1 {l\left( {\varphi \left( y \right)} \right){Y_n}\left( y \right)dy},\eqno(11)\]
or
\[{\lambda _n}{\varphi _n} = \int\limits_0^1 {K\left( y \right)l\left( \varphi  \right)\frac{{{Y_n}\left( y \right)}}{{K\left( y \right)}}dy} .\]
Consequently, ${\lambda _n}\,{\varphi _n}$ are the Fourier coefficients of the function $K\left( y \right)l\left( \varphi \right).$ Then from the Bessel inequality we get
\[\sum\limits_{n = 1}^\infty  {\lambda _n^2{{\left| {{\varphi _n}} \right|}^2}}  \le \int\limits_0^1 {K\left( y \right){{\left\{ {l\left( \varphi  \right)} \right\}}^2}dy}<\infty.\eqno(12)\]
From (7) and (12) the convergence of the series (10) and the fulfillment of condition (3) follows. Let us pass to the derivatives. In view of equation (8), we formally have
\[\left| {{}_CD_{0x}^\alpha u\left( {x,y} \right)} \right| \le \sum\limits_{n = 1}^\infty  {\left| {{}_CD_{0x}^\alpha {X_n}\left( x \right)} \right|\left| {{Y_n}\left( y \right)} \right|}  \le M{x^\beta }\sum\limits_{n = 1}^\infty  {\left| {{\lambda _n}{\varphi _n}} \right|\left| {{Y_n}\left( y \right)} \right|}  \le \]
\[ \le M{x^\beta }\sqrt {\sum\limits_{n = 1}^\infty  {{{\left( {\lambda _n^2{\varphi _n}} \right)}^2}} } \sqrt {\sum\limits_{n = 1}^\infty  {{{\left( {\frac{{{Y_n}\left( y \right)}}{{{\lambda _n}}}} \right)}^2}} },x>0.\eqno(13)\]
Further from (11) we have
\[{\varphi _n} = {\left( {\frac{1}{{{\lambda _n}}}} \right)^2}\int\limits_0^1 {K\left( y \right)l\left( {\varphi \left( y \right)} \right)l\left( {{Y_n}\left( y \right)} \right)dy}  = {\left( {\frac{1}{{{\lambda _n}}}} \right)^2}\int\limits_0^1 {l\left( {K\left( y \right)l\left( {\varphi \left( y \right)} \right)} \right){Y_n}\left( y \right)dy} ,\]
now we apply the Bessel inequality
\[\sum\limits_{n = 1}^\infty  {{{\left( {\lambda _n^2{\varphi _n}} \right)}^2}}  \le \int\limits_0^1 {K\left( y \right){{\left\{ {l\left( {K\left( y \right)l\left( \varphi  \right)} \right)} \right\}}^2}dy < \infty .}\eqno(14) \]
From (7) and (14) convergence (13) follows. So the following functional series
\[D_{0x}^\alpha u\left( {x,y} \right) = \sum\limits_{n = 1}^\infty  {D_{0x}^\alpha {X_n}\left( x \right){Y_n}\left( y \right)} ,\]
converges evenly. Uniform convergence of the next series
\[l\left( {u\left( {x,y} \right)} \right) = \sum\limits_{n = 1}^\infty  {{X_n}\left( x \right)l\left( {{Y_n}\left( y \right)} \right)}  = \frac{1}{{K\left( y \right)}}\sum\limits_{n = 1}^\infty  {{\lambda _n}{X_n}\left( x \right){Y_n}\left( y \right)} ,\]
shown in the same way.

\textbf{Theorem 2 is proved.}

\section*{3. Uniqueness of the solution }

Let us proceed to the proof of the uniqueness of the solution. The following theorem holds.

\textbf{Theorem 3.} If there is a solution to Problem A, then it is unique.

\textbf{ Proof.} Let the function $u\left( {x,y} \right)$ be a solution to Problem A with zero initial and boundary conditions. Consider its Fourier coefficients in terms of the system of eigenfunctions of problem (4)
\[{u_n}\left( x \right) = \int\limits_0^1 {\frac{{u\left( {x,y} \right){Y_n}(y)dy}}{{K\left( y \right)}}} ,{\mkern 1mu} {\kern 1pt} n = 1,2,...,\]
it is not hard to show that ${u_n}\left( x \right)$ is a solution to the problem
\[\left\{ {\begin{array}{*{20}{l}}
{{}_CD_{0x}^\alpha {u_n}\left( x \right) =  - {\lambda _n}{x^\beta }{u_n}\left( x \right),{\lambda _n} > 0,}\\
{{u_n}\left( 0 \right) = 0.}
\end{array}} \right.\]
This problem has only a trivial solution [14], i.e.,
\[\int\limits_0^1 {\frac{{u\left( {x,y} \right){Y_n}(y)dy}}{{K\left( y \right)}}}  = 0,n = 1,2,....\]
Because $\bar G\left( {y,\xi } \right)  $ is a symmetric, continuous function and there are integrals:
\[\int\limits_0^1 {{{\bar G}^2}\left( {y,\xi } \right)d\xi }  < \infty ,{\mkern 1mu} {\kern 1pt} \int\limits_0^1 {{{\bar G}^2}\left( {y,\xi } \right)dy}  < \infty ,{\mkern 1mu} {\kern 1pt} \int\limits_0^1 {\int\limits_0^1 {{{\bar G}^2}\left( {y,\xi } \right)dyd\xi } }  < \infty ,{\mkern 1mu} {\kern 1pt} {\lambda _n} > 0,{\mkern 1mu} {\kern 1pt} \forall n,\]
then it follows from Mercer's theorem that
\[\bar G\left( {y,\xi } \right) = \sum\limits_{n = 1}^\infty  {\frac{{\overline {{Y_n}} \left( y \right)\overline {{Y_n}} \left( \xi  \right)}}{{{\lambda _n}}}} .\]	
Then we have
\[\frac{{u\left( {x,y} \right)}}{{\sqrt {K\left( y \right)} }} = \int\limits_0^1 {\bar G\left( {y,\xi } \right)\left( {\sqrt {K\left( \xi  \right)} l\left( {u\left( {x,\xi } \right)} \right)} \right)d\xi }  = \]
\[ = \int\limits_0^1 {\sum\limits_{n = 1}^\infty  {\frac{{\overline {{Y_n}} \left( y \right)\overline {{Y_n}} \left( \xi  \right)}}{{{\lambda _n}}}} \left( {\sqrt {K\left( \xi  \right)} l\left( {u\left( {x,\xi } \right)} \right)} \right)d\xi }  = \]
\[ = \sum\limits_{n = 1}^\infty  {\frac{{{Y_n}\left( y \right)}}{{{\lambda _n}\sqrt {K\left( y \right)} }}} \int\limits_0^1 {\frac{{{Y_n}\left( \xi  \right)\sqrt {K\left( \xi  \right)} }}{{\sqrt {K\left( \xi  \right)} }}l\left( {u\left( {x,\xi } \right)} \right)d\xi }  = \]
since the series converges uniformly, the signs of integration and the sum can be rearranged
\[ = \sum\limits_{n = 1}^\infty  {\frac{{{Y_n}\left( y \right)}}{{{\lambda _n}\sqrt {K\left( y \right)} }}} \int\limits_0^1 {{Y_n}\left( \xi  \right)l\left( {u\left( {x,\xi } \right)} \right)d\xi }  = \]
\[ = \sum\limits_{n = 1}^\infty  {\frac{{{Y_n}\left( y \right)}}{{{\lambda _n}\sqrt {K\left( y \right)} }}} \int\limits_0^1 {l\left( {{Y_n}\left( \xi  \right)} \right)u\left( {x,\xi } \right)d\xi }  = \]
\[ = \sum\limits_{n = 1}^\infty  {\frac{{{Y_n}\left( y \right)}}{{\sqrt {K\left( y \right)} }}} \int\limits_0^1 {\frac{{{Y_n}\left( \xi  \right)u\left( {x,\xi } \right)}}{{K\left( \xi  \right)}}d\xi } =0, \]
from here
\[u\left( {x,y} \right) \equiv 0.\]
	
\textbf{Theorem 3 is proved.}

\begin{center}
\textbf{References}
\end{center}
1. A. N. Bogolyubov, A. A. Koblikov, D. D. Smirnova, and N. E. Shapkina. Mathematical modelling of media
with time dispersion using fractional derivatives. Matem. Mod.25(12), 50–64 (2013)\\
2. A. M. Nakhushev.Fractional Calculus and Its Application(Fizmatlit, Moscow, 2003) [in Russian]\\
3. A.N. Artyushin.  Fractional integral inequalities and their applications to degenerate differential equations with the Сaputo fractional derivative. Siberian Mathematical Journal , 62 (2), 2020. pp. 208 - 221.\\
4. A. A. Kilbas, Megumi Saigo. Solution in closed form of a class of linear differential equations of fractional order.
Different. Equat., 33, № 2, 194 – 204 (1997).\\
5. Samko S. G., Kilbas A. A., Marichev O. I. Fractional integrals and derivatives and some of their applications. - Minsk. Science and technology. 1987. - 688 p.[in Russian]\\
6. Oldham К. В., Spanier J. The fractional calculus. New York; London. 1974.\\
7. Wiener K.  Wiss. Z. Univ. Halle Math. Natur. Wiss. R.  1983. 32 (1), 1983. pp. 41 - 46.\\
8. Karimov E., Ruzhansky M., Toshtemirov B.  Solvability of the boundary-value problem for a mixed equation involving hyper-Bessel fractional differential operator and bi-ordinal Hilfer fractional derivative. Mathematical Methods in the Applied Sciences. 41(1), 2023, pp. 54-77.\\
9. Ruzhansky M., Torebek B. T., Turmetov B.  Well-posedness of Tricomi-Gellerstedt-Keldysh-type fractional elliptic problems. Journal of Integral Equations and Applications. 34 (3), 2022, pp. 373 - 387.\\
10. Smadiyeva A.G. Initial-boundary value problem for the time-fractional degenerate diffusion equation. Journal of Mathematics, Mechanics and Computer Science, [S.l.], v. 113, n. 1, mar. 2022. ISSN 2617-4871.
 DOI: https: // doi.org /10.26577 /JMMCS. 2022. v113.i1.04 \\
11. Smadiyeva A.G.   Well-posedness of the initial-boundary value problems for the
time-fractional degenerate diffusion equations.  Bulletin of the Karaganda University. Mathematics series.107(3),2022,pp.145-151.DOI: 10.31489/2022M3/145-151\\
12. Urinov A. K., Usmonov D. A. Non-local initial-boundary value problem for a degenerate fourth-order equation with a fractional Gerasimov-Caputo derivative. Vestnik KRAUNC. Fiz.-mat. nauki. 2023, 42: 1, 123-139. EDN: INZPHJ. https://doi.org/10.26117/2079-6641-2023-42-1-123-139.\\
13. Irgashev, B.Y. A Boundary Value Problem with Conjugation Conditions for a Degenerate Equation with the Caputo Fractional Derivative. Russ Math. 66, 24–31 (2022). https://doi.org/10.3103/S1066369X2204003X\\
14. Kilbas, Anatoly A.; Srivastava, Hari M.; Trujillo, Juan J. Theory and applications of fractional differential equations. North-Holland Mathematics Studies, 204. Elsevier Science B.V., Amsterdam, 2006.\\
15.Ravshan Ashurov, Yusuf Fayziev, Muattar Khudoykulova. Forward and Inverse Problems for Subdiffusion Equation with Time-Dependent Coefficients. 29 March, 2023. arXiv:2304.00998\\
16. Irgashev B.Y. Initial-boundary problem for degenerate high order equation with fractional derivative. Indian J Pure Appl Math 53, 170–180 (2022). https://doi.org/10.1007/s13226-021-00088-7\\
17.Irgashev B.Y. Mixed problem for higher-order equations with fractional derivative and degeneration in both variables. Ukrains’ kyi Matematychnyi Zhurnal. 2022. Vol.74. N 10. P. 1328-1338\\
18. Lotfi Boudabsa, Thomas Simon. Some Properties of the Kilbas-Saigo Function. Mathematics. 217 (9), 2021. https://doi.org/10.3390/math9030217
 \end{document}